\pdfoutput=1
\RequirePackage{ifpdf}
\ifpdf 
\documentclass[pdftex]{sigma}
\else
\documentclass{sigma}
\fi

\usepackage{stmaryrd}

\newtheorem*{remarks}{Remarks}
\newtheorem*{Remark}{Remark}

\begin{document}

\newcommand{\arXivNumber}{1305.6946}

\allowdisplaybreaks

\renewcommand{\PaperNumber}{072}

\FirstPageHeading

\ShortArticleName{The GraviGUT Algebra Is not a~Subalgebra of $E_8$}

\ArticleName{The GraviGUT Algebra Is not a~Subalgebra of $\boldsymbol{E_8}$,\\
but $\boldsymbol{E_8}$ Does Contain an Extended GraviGUT Algebra}

\Author{Andrew DOUGLAS~$^\dag$ and Joe REPKA~$^\ddag$}

\AuthorNameForHeading{A.~Douglas and J.~Repka}

\Address{$^\dag$~CUNY Graduate Center and New York City College of Technology,\\
\hphantom{$^\dag$}~City University of New York, USA}
\EmailD{\href{mailto:adouglas2@gc.cuny.edu}{adouglas2@gc.cuny.edu}}

\Address{$^\ddag$~Department of Mathematics, University of Toronto, Canada}
\EmailD{\href{mailto:repka@math.toronto.edu}{repka@math.toronto.edu}}

\ArticleDates{Received April 04, 2014, in f\/inal form July 03, 2014; Published online July 08, 2014}

\Abstract{The (real) GraviGUT algebra is an extension of the $\mathfrak{spin}(11,3)$ algebra by a~$64$-dimensional Lie
algebra, but there is some ambiguity in the literature about its def\/inition.
Recently, Lisi constructed an embedding of the GraviGUT algebra into the quaternionic real form of~$E_8$.
We clarify the def\/inition, showing that there is only one possibility, and then prove that the GraviGUT algebra cannot
be embedded into any real form of~$E_8$.
We then modify Lisi's construction to create true Lie algebra embeddings of the extended GraviGUT algebra into~$E_8$.
We classify these embeddings up to inner automorphism.}


\Classification{17B05; 17B10; 17B20; 17B25; 17B81}

\section{Introduction}

The Standard Model of particle physics, with gauge group $\text{U}(1) \times \text{SU}(2) \times \text{SU}(3)$, attempts
to describe all particles and all forces, except gravity.
Grand Unif\/ied Theories (GUT) attempt to unify the forces and particles of the Standard Model.
The three main GUTs are Georgi and Glashow's $\text{SU}(5)$ theory, Georgi's $\text{Spin}(10)$ theory, and the
Pati--Salam model based on the Lie group $\text{SU}(2) \times \text{SU}(2) \times \text{SU}(4)$~\cite{bh}.

In~\cite{lisi}, Lisi attempts to construct a~unif\/ication which includes gravity.
In this construction, Lisi f\/irst embeds gravity and the standard model into $\mathfrak{spin}(11,3)$.
He then embeds $\mathfrak{spin}(11,3)$ together with the positive chirality $64$-dimensional $\mathfrak{spin}(11,3)$
irrep into the quaternionic real form of~$E_8$.
Lisi refers to the embedded Lie algebra as the GraviGUT algebra.
For a~description of Lisi's theory see~\cite{lisi} or~\cite{lisi2}.
For a~critique of Lisi's theory involving the GraviGUT algebra see~\cite{sgd}.
We note that the GraviGUT algebra was f\/irst introduced by Nesti and Percacci~\cite{nesti}.

In Section~\ref{Section4} of~\cite{mkrtchyan}, it is observed that one of the Lie algebras associated with a~point in Vogel's plane
(cf.~\cite{vogel}) has the same dimension as the GraviGUT algebra, and it is hypothesized there that these two
algebras are in fact isomorphic.
It would certainly be interesting to identify the algebra in question.

Because the exposition of~\cite{lisi} ref\/lects the process of exploring the possible realizations of the GraviGUT
algebra, there is some potential for confusion about its def\/inition.
Lisi f\/irst describes the algebraic structure of $\mathfrak{spin}(11,3)$ and the action of $\mathfrak{spin}(11,3)$ on its
$64$-dimensional irrep~$V$ in equations $(3.9)$ and $(3.10)$ from \cite{lisi}, respectively.
He then notes that the structure of the GraviGUT algebra {\it could} be completed by def\/ining a~trivial Lie bracket on~$V$.

However, the actual algebraic structure of his explicit realization of~$V$
is not established until equations $(4.3)$ and $(4.4)$ from~\cite{lisi},
when Lisi describes his embedding; it is here that we see that~$V$ is not abelian relative to the Lie
bracket inherited from~$E_{8}$.
In fact, with this def\/inition of the Lie bracket, the subspace of $E_{8}$ spanned by the embedded copies of
$\mathfrak{spin}(11,3)$ and~$V$ is not a~Lie algebra at all: it is not closed under the bracket.
Lisi acknowledges this in his remark near the end of \S~3 that ``The word `algebra' is used here in a~generalized sense''.

Theorem~\ref{pijjk} of the present paper actually shows that the {\it only} way to extend the usual bracket on $
\mathfrak{so}(14)_{\mathbb{C}} $ and its action on~$V$ to make $ \mathfrak{so}(14)_{\mathbb{C}} \inplus V$ into a~Lie
algebra is to require~$V$ to be abelian.
In particular, the only possible def\/inition of the (complexif\/ied) GraviGUT algebra as a~Lie algebra
is $\mathfrak{so}(14)_{\mathbb{C}} \inplus V$, where~$V$ is a~$64$-dimensional abelian ideal which is irreducible under the
action of $ \mathfrak{so}(14)_{\mathbb{C}}$.

Once the structure of the complexif\/ied GraviGUT algebra is specif\/ied, it is not dif\/f\/icult to show that it cannot be
embedded into the complex algebra $E_8$; cf.\ Corollary~\ref{Cor:OnlyAbelian}.
Hence, the (real) GraviGUT algebra cannot be embedded into the quaternionic real form of $E_8$, or any other real form
of~$E_8$.

However, the operators in $E_{8}$ described by Lisi do generate a~larger Lie algebra, which contains an additional
$14$-dimensional ideal.
We call this larger algebra the extended GraviGUT algebra.

We modify Lisi's construction to create true Lie algebra embeddings of the extended Gravi\-GUT algebra into~$E_8$.
The (complex) \emph{extended GraviGUT algebra} is a~nonabelian, nilpotent extension of $ \mathfrak{so}(14)_{\mathbb{C}}
$ by a~$78$-dimensional $ \mathfrak{so}(14)_{\mathbb{C}} $-representation.
This $78$-dimensional representation is composed of a~$64$-dimensional irrep and the standard $14$-dimensional $
\mathfrak{so}(14)_{\mathbb{C}} $-irrep.
Its precise structure is described in Section~\ref{Section6}, but we do note here that the (complexif\/ied) Gravi\-GUT algebra is
a~quotient of the extended Gravi\-GUT algebra.
We classify these embeddings up to inner automorphism.

The article is organized as follows.
Section~\ref{cb} contains relevant background on Lie algebras and their representations: in particular, it deals with the
complex, simple Lie algebras $\mathfrak{so}(14)_{\mathbb{C}}$ and~$E_8$.
Section~\ref{Section3} presents additional notation and terminology.
In Section~\ref{Section4} we describe the classif\/ication of embeddings of $ \mathfrak{so}(14)_{\mathbb{C}} $ into $E_8$, which will
be used in the following section.
In Section~\ref{last} we determine the only possible def\/inition of the GraviGUT algebra and also establish that the
complexif\/ied GraviGUT algebra cannot be embedded into the complex algebra~$E_8$.
Finally, in Section~\ref{Section6} we classify the embeddings of the extended GraviGUT algebra into~$E_8$.

\section[The complex Lie algebras $ \mathfrak{so}(14)_{\mathbb{C}} $ and $E_8$,\\and their representations]
{The complex Lie algebras $\boldsymbol{\mathfrak{so}(14)_{\mathbb{C}}}$ and $\boldsymbol{E_8}$,\\and their representations}
\label{cb}

The special orthogonal algebra $ \mathfrak{so}(14)_{\mathbb{C}} $ is the complexif\/ication of $\mathfrak{spin}(11,3)$.
It is the Lie algebra of complex $14\times 14$ matrices~$N$ satisfying $N^{\rm tr}=-N$.
The dimension of $ \mathfrak{so}(14)_{\mathbb{C}} $ is $91$ and its rank is~$7$.
The Lie group corresponding to $ \mathfrak{so}(14)_{\mathbb{C}} $ arises naturally as the symmetry group of a~projective
space over $\mathbb{R}$~\cite{baez}.

$E_8$ is the complex, exceptional Lie algebra of rank~$8$.
It is $248$-dimensional.
Like $ \mathfrak{so}(14)_{\mathbb{C}} $, $E_8$ has a~close connection to the Riemannian geometry of projective spaces
(for details, we refer the reader to~\cite{baez}).

Let $\mathfrak{g}$ denote $ \mathfrak{so}(14)_{\mathbb{C}} $ or~$E_8$.
Let $k=7$ or $8$ when $ \mathfrak{g} = \mathfrak{so}(14)_{\mathbb{C}} $ or $E_8$, respectively.
We may def\/ine $\mathfrak{g}$ by a~set of generators $\{H_i, X_i, Y_i\}_{1\leq i \leq k}$ together with the
Chevalley--Serre relations~\cite{humphreys}:
\begin{gather*}
 [H_i,H_j]=0, \qquad  [H_i,X_j]= M^{\mathfrak{g}}_{ji} X_j,
\\
 [H_i,Y_j]= -M^{\mathfrak{g}}_{ji} Y_j, \qquad  [X_i,Y_j]= \delta_{ij} H_i,
\\
 (\text{ad}\,X_i)^{1-M^{\mathfrak{g}}_{ji}}(X_{j}) =0,
\qquad
(\text{ad}\, Y_i)^{1-M^{\mathfrak{g}}_{ji}}(Y_j)=0,
\quad
\text{when}
\quad
i\neq j.
\end{gather*}
Here $1 \leq i,j \leq k$, and $M^{\mathfrak{g}}$ is the Cartan matrix of $\mathfrak{g}$.
The $X_i$, for $1\leq i \leq k$, correspond to the simple roots.
We write~$H$ for the Cartan subalgebra spanned by $\{H_{i}\}$.

\begin{figure}[h!]\centering
\includegraphics{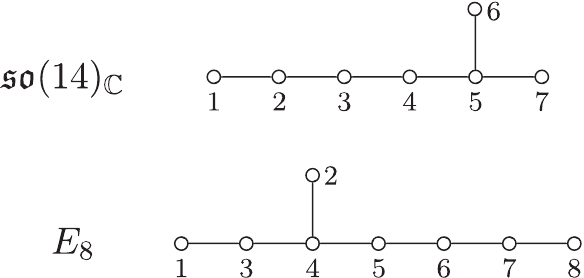}
\caption{Dynkin diagrams of $ \mathfrak{so}(14)_{\mathbb{C}} $ and~$E_8$.}
\label{ddd}
\end{figure}

For future reference, the Dynkin diagrams of $ \mathfrak{so}(14)_{\mathbb{C}} $ and $E_8$, indicating the numbering of
simple roots, are given in Fig.~\ref{ddd}.

We now brief\/ly describe the f\/inite-dimensional, irreducible representations (irreps) of $ \mathfrak{so}(14)_{\mathbb{C}}
$ and $E_8$, with $\mathfrak{g}$ and~$k$ def\/ined as above.
For $i=1,\ldots,k$, def\/ine $\alpha_i,\lambda_i \in H^*$ by $\alpha_i(H_j)= M^{\mathfrak{g}}_{ji}$, and
$\lambda_i(H_j)=\delta_{ij}$, where $M^{\mathfrak{g}}$ is the Cartan matrix of $\mathfrak{g}$.
The $\lambda_i$ are the {\it fundamental weights}, and their indexing corresponds with that of the Dynkin diagram of
type $ \mathfrak{so}(14)_{\mathbb{C}} $ or $E_8$ in Fig.~\ref{ddd}.

For each $\lambda= m_1 \lambda_1+\dots+m_k \lambda_k \in H^*$ with nonnegative integers $m_1,\ldots,m_k$, there exists an
irrep of $\mathfrak{g}$ with highest weight~$\lambda$, denoted $V_{\mathfrak{g}}(\lambda)$.
The irreps $V_{\mathfrak{g}}(\lambda_i)$ for $1\leq i \leq k$ are the {\it fundamental representations}.
Each irrep of $\mathfrak{g}$ is equivalent to $V_{\mathfrak{g}}(\lambda)$, where $\lambda=m_1\lambda_1+\dots+m_k\lambda_k$
for some nonnegative integers $m_1,\ldots, m_k$.

\section{Additional def\/initions and notation}
\label{Section3}

The following def\/initions and notation will be used in this article:
\begin{itemize}\itemsep=0pt
\item For $ 1\leq a_i \leq 8$, let $X_{a_i}$ correspond to the $a_i$th simple root of~$E_8$.
We then def\/ine
\begin{gather*}
X_{a_1,a_2,a_3,\ldots,a_ m} \equiv [[\ldots[[X_{a_{1}},X_{a_{2}}],X_{a_{3}}],\ldots],X_{a_{m}}].
\end{gather*}
$Y_{a_1,a_2,a_3,\ldots,a_m}$ is def\/ined analogously.

\item Let $\varphi: \mathfrak{so}(14)_{\mathbb{C}} \hookrightarrow E_{8}$ be an embedding.
Further, let~$W$ be an element of $E_{8}$.
Then, $[W]_{\varphi(\mathfrak{so}(14)_{\mathbb{C}})}$ is the $ \mathfrak{so}(14)_{\mathbb{C}} $ representation generated
by~$W$ with respect to the adjoint action of $\varphi(\mathfrak{so}(14)_{\mathbb{C}})$.
When the embedding~$\varphi$ is clear, as will be the case below, we simply write
$[W]_{\mathfrak{so}(14)_{\mathbb{C}}}$.

\item Let~$\varphi$ and~$\varrho$ be Lie algebra embeddings of $\mathfrak{g}'$ into $\mathfrak{g}$.
Then~$\varphi$ and~$\varrho$ are {\it equivalent} if there is an inner automorphism $\rho: \mathfrak{g} \rightarrow
\mathfrak{g}$ such that $\varphi= \rho \circ \varrho$, and we write
\begin{gather*}
\varphi \sim \varrho.
\end{gather*}
\item Two embeddings~$\varphi$ and~$\varrho$ of $\mathfrak{g}'$ into $\mathfrak{g}_{}$ are \emph{linearly equivalent} if
for each representation $\pi: \mathfrak{g} \rightarrow \mathfrak{gl}(V)$ the induced $\mathfrak{g}'$ representations
$\pi \circ \varphi$, $\pi \circ \varrho$ are equivalent, and we write
\begin{gather*}
\varphi \sim_L \varrho.
\end{gather*}
Clearly equivalence of embeddings implies linear equivalence, but the converse is not in general true.
\end{itemize}

We def\/ine equivalence and linear equivalence of subalgebras much as we did for embeddings:
\begin{itemize}\itemsep=0pt

\item Two subalgebras $\mathfrak{g}'$ and $\mathfrak{g}''$ of $\mathfrak{g}$ are equivalent if there is an inner
automorphism~$\rho$ of $\mathfrak{g}$ such that $\rho(\mathfrak{g}')=\mathfrak{g}''$.
\item Two subalgebras $\mathfrak{g}'$ and $\mathfrak{g}''$ of $\mathfrak{g}$ are linearly equivalent if for every
representation $\pi:\mathfrak{g} \rightarrow \mathfrak{gl}(V)$ the subalgebras $ \pi (\mathfrak{g}')$, $ \pi
(\mathfrak{g}'')$ of $\mathfrak{gl}(V)$ are conjugate under ${\rm GL}(V)$.
\end{itemize}

\section[Embedding $ \mathfrak{so}(14)_{\mathbb{C}} $ into $E_{8}$]
{Embedding $\boldsymbol{\mathfrak{so}(14)_{\mathbb{C}}}$ into $\boldsymbol{E_{8}}$}
\label{Section4}

In~\cite{dkr}, the authors presented the following well-known ``natural'' embedding of $ \mathfrak{so}(14)_{\mathbb{C}}$
into $E_{8}$:
\begin{gather}
\varphi:   \ \mathfrak{so}(14)_{\mathbb{C}}  \hookrightarrow  E_{8},
\nonumber
\\
\phantom{\varphi:{}} \  H_{8-i}   \mapsto   H_{i+1},
\nonumber
\\
\phantom{\varphi:{}}  \ X_{8-i}   \mapsto   X_{i+1},
\nonumber
\\
\phantom{\varphi:{}}  \ Y_{8-i}   \mapsto   Y_{i+1},
\label{embb}
\end{gather}
where $1\leq i\leq 7$.
This embedding may be visualized as a~``natural'' subgraph of the Dynkin diagram of $E_{8}$ which is isomorphic to the
Dynkin diagram of $ \mathfrak{so}(14)_{\mathbb{C}} $ (see Fig.~\ref{ddd}).

In~\cite{min}, Minchenko showed that there is a~unique subalgebra isomorphic to $ \mathfrak{so}(14)_{\mathbb{C}} $ in
$E_8$, up to inner automorphism.
Hence, the only way to get new embeddings of $ \mathfrak{so}(14)_{\mathbb{C}} $ into $E_8$ other than the~$\varphi$
described in equation~\eqref{embb} is to compose~$\varphi$ with an outer automorphism of $ \mathfrak{so}(14)_{\mathbb{C}} $.
However, it was shown by the authors in~\cite{dkr} that outer automorphisms of $ \mathfrak{so}(14)_{\mathbb{C}} $ do not
produce new embeddings of $ \mathfrak{so}(14)_{\mathbb{C}} $ into~$E_8$.
We thus have the following theorem~\cite{dkr}.

\begin{theorem}
\label{cla}
The map $\varphi: \mathfrak{so}(14)_{\mathbb{C}} \hookrightarrow E_8$ defined in equation~\eqref{embb} is the unique
embedding of $ \mathfrak{so}(14)_{\mathbb{C}} $ into $E_8$, up to inner automorphism.
\end{theorem}

\section[The GraviGUT algebra is not a~subalgebra of $E_8$]{The GraviGUT algebra is not a~subalgebra of $\boldsymbol{E_8}$}
\label{last}

\begin{theorem}\label{pijjk}
Consider a~sum of complex vector spaces $ \mathfrak{so}(14)_{\mathbb{C}} \oplus V$, where~$V$ is a~$64$-dimensional
space.
Suppose a~Lie bracket is defined which gives the usual structure to the $ \mathfrak{so}(14)_{\mathbb{C}} $ subspace and
in such a~way that the brackets $[X,v]$, for $X \in \mathfrak{so}(14)_{\mathbb{C}} $ and $v \in V$, define an action of~$ \mathfrak{so}(14)_{\mathbb{C}} $ on~$V$ under which~$V$ becomes an irreducible $ \mathfrak{so}(14)_{\mathbb{C}}
$-module.
Then the {\it only} way to extend this bracket to make $ \mathfrak{so}(14)_{\mathbb{C}} \inplus V$ into a~Lie algebra is
to put the abelian structure on~$V$, i.e., $[v,v'] =0$, for all $v,v' \in V$.

In particular, the only possible definition of the $($complexified$)$ GraviGUT algebra as a~Lie algebra
is $\mathfrak{so}(14)_{\mathbb{C}} \inplus V$, where~$V$ is a~$64$-dimensional abelian ideal which is irreducible under the
action of $ \mathfrak{so}(14)_{\mathbb{C}} $.
\end{theorem}

\begin{proof}
Let~$V$ be a~$64$-dimensional $ \mathfrak{so}(14)_{\mathbb{C}} $-irrep.
Then $V\cong V(\lambda_6)$ or $V\cong V(\lambda_7)$.
Consider the tensor product decompositions:
\begin{gather*}
V(\lambda_6) \otimes V(\lambda_6) \cong V(2\lambda_6) \oplus V(\lambda_5)\oplus V(\lambda_3)\oplus V(\lambda_1),
\\
V(\lambda_7) \otimes V(\lambda_7) \cong V(2\lambda_7) \oplus V(\lambda_5)\oplus V(\lambda_3)\oplus V(\lambda_1).
\end{gather*}
Since $V(\lambda_6)$ does not occur in the tensor product decomposition of $V(\lambda_6) \otimes V(\lambda_6)$,
$V(\lambda_7)$ does not occur in the tensor product decomposition of $V(\lambda_7) \otimes V(\lambda_7)$, and neither
decomposition contains a~$91$-dimensional irrep, we cannot have a~nontrivial product $V \otimes V \rightarrow V$ or $V
\otimes V \rightarrow \mathfrak{so}(14)_{\mathbb{C}} $.
Hence, we cannot have a~nontrivial product $V \times V \rightarrow V$ or $V \times V \rightarrow
\mathfrak{so}(14)_{\mathbb{C}} $ or $V \times V \rightarrow (\mathfrak{so}(14)_{\mathbb{C}} \inplus V)$.
The only possible def\/inition of a~Lie algebra structure is to make~$V$ an abelian subalgebra.
\end{proof}

\begin{corollary}
\label{Cor:OnlyAbelian}
The GraviGUT algebra cannot be embedded into the quaternionic real form of $E_8$, or any other real form of~$E_8$.
\end{corollary}
\begin{proof}
The maximal dimension of an abelian subalgebra of $E_8$ is $36$~\cite{mal}.
This implies that~$V$ cannot be a~subalgebra of $E_8$, and that the complexif\/ied GraviGUT algebra cannot be embedded
into~$E_8$.
The result follows.
\end{proof}

\begin{remarks}\rm \quad

\begin{enumerate}\itemsep=0pt
\item[1.] We also note that if we have a~$64$-dimensional representation~$V$ that is not irreducible, then embedding $
\mathfrak{so}(14)_{\mathbb{C}} \inplus V$ into $E_8$ is still not possible.
The summands in the direct sum decomposition of $E_8$ as an $\varphi(\mathfrak{so}(14)_{\mathbb{C}})$-module have
dimensions $1$, $14$, $14$, $64$, $64$, $91$, as we will see below in equation~\eqref{hhhhkj}.
Hence, the only $64$-dimensional $ \mathfrak{so}(14)_{\mathbb{C}} $ submodules of $E_8$ are irreducible.

\item[2.] The subspace of $E_{8}$ def\/ined by Lisi in his equations $(4.3)$, $(4.4)$ (see~\cite{lisi})
is not closed under the Lie bracket it inherits from~$E_{8}$.
\end{enumerate}
\end{remarks}

\section[The extended GraviGUT algebra in $E_8$]{The extended GraviGUT algebra in $\boldsymbol{E_8}$}
\label{Section6}

In~\cite{dkr}, the authors computed the following decomposition of $E_8$ with respect to the adjoint action of
$\varphi(\mathfrak{so}(14)_{\mathbb{C}})$:
\begin{gather}
E_8\cong_{\mathfrak{so}(14)_{\mathbb{C}}}   V_{}(\lambda_2)  \oplus  V_{}(\lambda_1)  \oplus  V_{}(\lambda_6) \oplus
V_{}(\lambda_7)  \oplus  V_{}(\lambda_1) \oplus  V_{}(0)
\label{hhhhkj}
\\
\phantom{E_8}
\cong_{\mathfrak{so}(14)_{\mathbb{C}}}  [X_{74}]_{\mathfrak{so}(14)_{\mathbb{C}}}\! \oplus
[X_{120}]_{\mathfrak{so}(14)_{\mathbb{C}}} \! \oplus  [Y_{1}]_{\mathfrak{so}(14)_{\mathbb{C}}}\! \oplus [X_{112}]_{\mathfrak{so}(14)_{\mathbb{C}}}\!
\oplus  [Y_{97}]_{\mathfrak{so}(14)_{\mathbb{C}}}\! \oplus [H]_{\mathfrak{so}(14)_{\mathbb{C}}},
\nonumber
\end{gather}
where
\begin{gather*}
  X_{74} = X_{4,5,6,7,8,2,3,4,5,6,7},
\\
  X_{112}  = -X_{3,4,2,1,5,4,3,6,5,4,7,2,6,5,8,7,6,4,5,3,4,2},
\\
  X_{120} =  X_{8,7,6,5,4,3,2,1,4,5,6,7,3,4,5,6,2,4,5,3,4,2,1,3, 4,5,6,7,8},
\\
  Y_{97} =  -Y_{5,4,2,3,6,4,1,3,5,4,7,2,6,5,4,3,1},
\\
 H =  4H_1+5H_2+7H_3+10H_4+
8H_5+6H_6+4H_7+2H_8.
\end{gather*}

\begin{lemma}
\label{l1}
The following are $78$-dimensional, nonabelian nilpotent subalgebras of $E_8$:
\begin{gather*}
[Y_{1}]_{\mathfrak{so}(14)_{\mathbb{C}}}\oplus [Y_{97}]_{\mathfrak{so}(14)_{\mathbb{C}}}, \qquad
[X_{112}]_{\mathfrak{so}(14)_{\mathbb{C}}}\oplus [X_{120}]_{\mathfrak{so}(14)_{\mathbb{C}}}.
\end{gather*}
Note that the sums are direct as $ \mathfrak{so}(14)_{\mathbb{C}} $ irreps, but not as subalgebras of~$E_8$.
Further, the subalgebras $[Y_{97}]_{\mathfrak{so}(14)_{\mathbb{C}}}$ and $[X_{120}]_{\mathfrak{so}(14)_{\mathbb{C}}}$ of
$E_8$ are abelian.
\end{lemma}
\begin{proof}
The authors showed in~\cite{dkr} that $[Y_{97}]_{\mathfrak{so}(14)_{\mathbb{C}}}$ and
$[X_{120}]_{\mathfrak{so}(14)_{\mathbb{C}}}$ are abelian subalgebras of~$E_8$.

The positive roots of $E_8$, as explicitly described in Appendix~\ref{roots}, give us a~triangular decomposition of
$E_8$: $E_{8,+} \oplus E_{8,0} \oplus E_{8,-}$.
In Appendix~\ref{roots} we also explicitly describe bases for the representations
$[X_{120}]_{\mathfrak{so}(14)_{\mathbb{C}}}$, $[X_{112}]_{\mathfrak{so}(14)_{\mathbb{C}}}$,
$[Y_{97}]_{\mathfrak{so}(14)_{\mathbb{C}}}$, and $[Y_1]_{\mathfrak{so}(14)_{\mathbb{C}}}$.
Since $\varphi(\mathfrak{so}(14)_{\mathbb{C}})$ $=$ $[X_{74}]_{\mathfrak{so}(14)_{\mathbb{C}}}$, we also give bases for
$\varphi(\mathfrak{so}(14)_{\mathbb{C}+})$ and $\varphi(\mathfrak{so}(14)_{\mathbb{C}-})$.
Each of these bases consists of all positive root vectors, or all negative root vectors.

The bases given in Appendix~\ref{roots} imply
\begin{gather*}
E_{8,+} = \varphi(\mathfrak{so}(14)_{\mathbb{C}+})  \oplus  [X_{112}]_{\mathfrak{so}(14)_{\mathbb{C}}}  \oplus
[X_{120}]_{\mathfrak{so}(14)_{\mathbb{C}}},
\\
E_{8,-} = \varphi(\mathfrak{so}(14)_{\mathbb{C}-})  \oplus  [Y_1]_{\mathfrak{so}(14)_{\mathbb{C}}}  \oplus
[Y_{97}]_{\mathfrak{so}(14)_{\mathbb{C}}}.
\end{gather*}
And of course
\begin{gather*}
  [[X_{112}]_{\mathfrak{so}(14)_{\mathbb{C}}} \oplus [X_{120}]_{\mathfrak{so}(14)_{\mathbb{C}}},
[X_{112}]_{\mathfrak{so}(14)_{\mathbb{C}}} \oplus [X_{120}]_{\mathfrak{so}(14)_{\mathbb{C}}}],
\\
  [[Y_1]_{\mathfrak{so}(14)_{\mathbb{C}}} \oplus [Y_{97}]_{\mathfrak{so}(14)_{\mathbb{C}}},
[Y_1]_{\mathfrak{so}(14)_{\mathbb{C}}} \oplus [Y_{97}]_{\mathfrak{so}(14)_{\mathbb{C}}}]
\end{gather*}
are subsets of $E_{8,+}$ and $E_{8,-}$, respectively.

In Appendix~\ref{roots}, we see that the positive root vector $X_{\alpha}$ is in the basis of
$[X_{112}]_{\mathfrak{so}(14)_{\mathbb{C}}}$ or $[X_{120}]_{\mathfrak{so}(14)_{\mathbb{C}}}$ if $\alpha^{1}\neq 0$,
where $\alpha^{1}$ is the f\/irst entry of~$\alpha$.
If $\alpha^{1}= 0$, then $X_{\alpha}$ is in the basis of $\varphi(\mathfrak{so}(14)_{\mathbb{C}+})$.

Thus, if $X_\alpha$ and $X_{\alpha'}$ are positive root vectors in the basis of
$[X_{112}]_{\mathfrak{so}(14)_{\mathbb{C}}}$ or $ [X_{120}]_{\mathfrak{so}(14)_{\mathbb{C}}}$ such that $[X_\alpha,
X_{\alpha'}]$ $\neq$~$0$, then this product is a~nonzero scalar multiple of $X_{\alpha+\alpha'}$, where $(\alpha
+\alpha')^1\neq 0$, so that $X_{\alpha+\alpha'}$ is an element of $[X_{112}]_{\mathfrak{so}(14)_{\mathbb{C}}}\oplus
[X_{120}]_{\mathfrak{so}(14)_{\mathbb{C}}}$.
Therefore,
\begin{gather*}
  [[X_{112}]_{\mathfrak{so}(14)_{\mathbb{C}}} \oplus [X_{120}]_{\mathfrak{so}(14)_{\mathbb{C}}},
[X_{112}]_{\mathfrak{so}(14)_{\mathbb{C}}} \oplus [X_{120}]_{\mathfrak{so}(14)_{\mathbb{C}}}]
\subseteq [X_{112}]_{\mathfrak{so}(14)_{\mathbb{C}}} \oplus [X_{120}]_{\mathfrak{so}(14)_{\mathbb{C}}}.
\end{gather*}

In a~similar manner we show
\begin{gather*}
  [[Y_1]_{\mathfrak{so}(14)_{\mathbb{C}}} \oplus [Y_{97}]_{\mathfrak{so}(14)_{\mathbb{C}}},
[Y_1]_{\mathfrak{so}(14)_{\mathbb{C}}} \oplus [Y_{97}]_{\mathfrak{so}(14)_{\mathbb{C}}}]
\subseteq [Y_1]_{\mathfrak{so}(14)_{\mathbb{C}}} \oplus [Y_{97}]_{\mathfrak{so}(14)_{\mathbb{C}}}.
\end{gather*}
Thus $[Y_{1}]_{\mathfrak{so}(14)_{\mathbb{C}}}\oplus [Y_{97}]_{\mathfrak{so}(14)_{\mathbb{C}}}$ and
$[X_{112}]_{\mathfrak{so}(14)_{\mathbb{C}}}\oplus [X_{120}]_{\mathfrak{so}(14)_{\mathbb{C}}}$ are subalgebras of~$E_8$.
Further, they are nilpotent since they are contained in $E_{8,-}$ or $E_{8,+}$, respectively.
\end{proof}

\begin{lemma}
\label{l2}
The following are not subalgebras of $E_8$:
\begin{gather*}
[Y_{97}]_{\mathfrak{so}(14)_{\mathbb{C}}}\oplus [X_{112}]_{\mathfrak{so}(14)_{\mathbb{C}}},
\qquad
[X_{120}]_{\mathfrak{so}(14)_{\mathbb{C}}}\oplus [Y_1]_{\mathfrak{so}(14)_{\mathbb{C}}}.
\end{gather*}
\end{lemma}
\begin{proof}
Referring to the bases of $[X_{120}]_{\mathfrak{so}(14)_{\mathbb{C}}}$ and $[Y_1]_{\mathfrak{so}(14)_{\mathbb{C}}}$
described in Appendix~\ref{roots}, we have $Y_{112}\in [Y_1]_{\mathfrak{so}(14)_{\mathbb{C}}}$, and of course
$X_{120}\in [X_{120}]_{\mathfrak{so}(14)_{\mathbb{C}}}$.
However, $[Y_{112}, X_{120}]$ is a~nonzero multiple of $X_{47}$, which is not in
$[X_{120}]_{\mathfrak{so}(14)_{\mathbb{C}}} \oplus [Y_1]_{\mathfrak{so}(14)_{\mathbb{C}}}$.
Hence $[X_{120}]_{\mathfrak{so}(14)_{\mathbb{C}}}\oplus [Y_1]_{\mathfrak{so}(14)_{\mathbb{C}}}$ is not a~subalgebra.
Similarly $[Y_{97}]_{\mathfrak{so}(14)_{\mathbb{C}}}\oplus [X_{112}]_{\mathfrak{so}(14)_{\mathbb{C}}}$ is not
a~subalgebra.
\end{proof}

We may now explicitly def\/ine the \emph{extended GraviGUT algebra} as follows.
As a~vector space, it is
\begin{gather*}
\mathfrak{so}(14)_{\mathbb{C}} \inplus (V(\lambda_6) \oplus V(\lambda_1)) \cong \mathfrak{so}(14)_{\mathbb{C}} \inplus
(V(\lambda_7) \oplus V(\lambda_1)).
\end{gather*}
The Lie algebra structure is inherited from that of $E_{8}$.
In particular, the following subalgebras are {\it not} direct sums as algebras, though the sums are direct as vector
spaces:
\begin{gather*}
V(\lambda_6)  +  V(\lambda_1)  \cong  [Y_1]_{\mathfrak{so}(14)_{\mathbb{C}}} +
[Y_{97}]_{\mathfrak{so}(14)_{\mathbb{C}}},
\\
V(\lambda_7)  +  V(\lambda_1)  \cong  [X_{112}]_{\mathfrak{so}(14)_{\mathbb{C}}} +
[X_{120}]_{\mathfrak{so}(14)_{\mathbb{C}}}.
\end{gather*}

We note further that not only are $ \mathfrak{so}(14)_{\mathbb{C}} \inplus (V(\lambda_6) + V(\lambda_1))$ and $
\mathfrak{so}(14)_{\mathbb{C}} \inplus (V(\lambda_7) + V(\lambda_1))$ isomorphic subalgebras, but they are equivalent
subalgebas of $E_8$, related by the Chevalley involution of~$E_8$.
Hence, we shall only consider $ \mathfrak{so}(14)_{\mathbb{C}} \inplus (V(\lambda_7) + V(\lambda_1))$.
It is signif\/icant to observe that the (complexif\/ied) GraviGUT algebra is a~quotient of the extended GraviGUT algebra.
The only distinction that can be made is that these two subalgebras of $E_{8}$ are {\it not} the same as $
\mathfrak{so}(14)_{\mathbb{C}} $-modules.

We now proceed to the classif\/ication of embeddings of the extended GraviGUT algebra into~$E_8$.
A~lift of $\varphi: \mathfrak{so}(14)_{\mathbb{C}} \hookrightarrow E_8$ to $ \mathfrak{so}(14)_{\mathbb{C}} \inplus
(V(\lambda_7) + V(\lambda_1))$ is completely determined by its def\/inition on highest weight vectors of $V(\lambda_7)$ and $V(\lambda_1)$.
Call these vectors~$u$ and~$v$, respectively.
Hence, for any $\alpha, \beta \in \mathbb{C}^*$, the following is a~lift of $\varphi: \mathfrak{so}(14)_{\mathbb{C}}
\hookrightarrow E_8$ to $ \mathfrak{so}(14)_{\mathbb{C}} \inplus (V(\lambda_7) + V(\lambda_1))$:
\begin{gather}
\widetilde{\varphi}_{\alpha,\beta}:  \ \mathfrak{so}(14)_{\mathbb{C}} \inplus (V(\lambda_7) + V(\lambda_1)) \hookrightarrow  E_8,
\nonumber
\\
\phantom{\widetilde{\varphi}_{\alpha,\beta}:{}} \
  u   \mapsto  \alpha X_{112},
\nonumber
\\
\phantom{\widetilde{\varphi}_{\alpha,\beta}:{}} \
  v   \mapsto  \beta X_{120}.
\label{ccv}
\end{gather}

\begin{theorem}
\label{t1}
All embeddings of the extended GraviGUT algebra into $E_8$ are given by $\widetilde{\varphi}_{\alpha,\beta}$, for all
$\alpha, \beta \in \mathbb{C}^*$.
These embeddings are classified according to the rule
\begin{gather*}
\widetilde{\varphi}_{\alpha,\beta}\sim \widetilde{\varphi}_{\alpha',\beta'} \ \Leftrightarrow \ \alpha'^2 \beta = \alpha^2\beta'.
\end{gather*}
\end{theorem}
\begin{proof}
First note that by Theorem~\ref{cla}, all embeddings of the extended GraviGUT algebras must come from lifts
of~$\varphi$, and hence, considering Lemmas~\ref{l1} and~\ref{l2}, equation~\eqref{ccv} def\/ines all embeddings of the
extended GraviGUT algebra into~$E_8$.

Def\/ine inner automorphisms of $E_8$ as follows:
\begin{alignat*}{4}
&  \rho: \quad &&   X_1  \mapsto  \alpha X_1, \qquad&&   Y_1  \mapsto   \frac{1}{\alpha}Y_1,  &
\\
&&&   X_i  \mapsto  X_i, \qquad &&  Y_i  \mapsto  Y_i, &
\\
&  \rho': \quad &&  X_1  \mapsto  \alpha' X_1, \qquad &&  Y_1  \mapsto    \frac{1}{\alpha'}Y_1,  &
\\
&&&   X_i  \mapsto  X_i, \qquad &&  Y_i  \mapsto  Y_i,&
\end{alignat*}
for $2\leq i \leq 8$.
Then $\widetilde{\varphi}_{\alpha,\beta}$ $=$ $\rho \circ \widetilde{\varphi}_{1, \frac{\beta}{\alpha^2}}$, and
$\widetilde{\varphi}_{\alpha',\beta'}$ $=$ $\rho' \circ \widetilde{\varphi}_{1, \frac{\beta'}{\alpha'^2}}$.
Hence $\widetilde{\varphi}_{\alpha,\beta}$ $\sim$ $\widetilde{\varphi}_{1, \frac{\beta}{\alpha^2}}$,
and $\widetilde{\varphi}_{\alpha',\beta'}$ $\sim$ $\widetilde{\varphi}_{1, \frac{\beta'}{\alpha'^2}}$.

If~$\vartheta$ is an inner automorphism of $E_8$ such that $\vartheta \circ \widetilde{\varphi}_{1,
\frac{\beta}{\alpha^2}}$ $=$ $\widetilde{\varphi}_{1, \frac{\beta'}{\alpha'^2}}$, then~$\vartheta$ f\/ixes $X_i$ and $Y_i$
for $2\leq i \leq 8$, and also $X_{112}$.
We have
\begin{gather*}
[\dots [X_{112}, Y_2],Y_4],Y_3],Y_5], Y_6], Y_7], Y_8], Y_4], Y_5],
\\
\qquad Y_2], Y_6], Y_4], Y_5], Y_3], Y_4], Y_7], Y_6], Y_5], Y_2], Y_4], Y_3]
=X_1,
\end{gather*}
so that~$\vartheta$ f\/ixes $X_1$.
Hence $\vartheta(X_{120})=X_{120}$, so that $\frac{\beta}{\alpha^2}=\frac{\beta'}{\alpha'^2}$.
The opposite implication is obvious.
Hence we have established
\begin{gather*}
\widetilde{\varphi}_{\alpha,\beta}\sim \widetilde{\varphi}_{\alpha,\beta} \Leftrightarrow \alpha'^2 \beta = \alpha^2\beta'.\tag*{\qed}
\end{gather*}
\renewcommand{\qed}{}
\end{proof}

\begin{Remark}\rm
Theorem~\ref{t1} implies, of course, that there are an inf\/inite number of embeddings of the extended GraviGUT algebra
into $E_8$, up to inner automorphism.
However, it is interesting to note that there is a~unique subalgebra of $E_8$ which is isomorphic to the extended
GraviGUT algebra up to inner automorphism.
\end{Remark}

\section{Conclusions}

In~\cite{lisi}, Lisi identif\/ied a~copy of $\mathfrak{spin}(11,3)$ in the quaternionic real form of $E_8$ and
a~$64$-dimensional subspace on which $\mathfrak{spin}(11,3)$ acts irreducibly.
His hope was to embed the GraviGUT algebra.
However, the subspace spanned by these spaces is not closed under the Lie bracket.

We proved that the only possible Lie algebra structure on the GraviGUT algebra has a~trivial (abelian) bracket on the
$64$-dimensional subspace.
In particular, the GraviGUT algebra cannot be embedded into any real form of~$E_8$.
We then modif\/ied Lisi's construction to create true Lie algebra embeddings of the extended GraviGUT algebra into~$E_8$.
We classif\/ied these embeddings up to inner automorphism.

\appendix

\section[\protect{The representations $[X_{74}]_{\mathfrak{so}(14)_{\mathbb{C}}}$, $[X_{120}]_{\mathfrak{so}(14)_{\mathbb{C}}}$,
$[X_{112}]_{\mathfrak{so}(14)_{\mathbb{C}}}$, $[Y_{97}]_{\mathfrak{so}(14)_{\mathbb{C}}}$, and $[Y_1]_{\mathfrak{so}(14)_{\mathbb{C}}}$}]
{The representations $\boldsymbol{[X_{74}]_{\mathfrak{so}(14)_{\mathbb{C}}}}$, $\boldsymbol{[X_{120}]_{\mathfrak{so}(14)_{\mathbb{C}}}}$,
$\boldsymbol{[X_{112}]_{\mathfrak{so}(14)_{\mathbb{C}}}}$,\\
$\boldsymbol{[Y_{97}]_{\mathfrak{so}(14)_{\mathbb{C}}}}$, and $\boldsymbol{[Y_1]_{\mathfrak{so}(14)_{\mathbb{C}}}}$}
\label{roots}

{\sloppy In this appendix
 we describe the representations $[X_{74}]_{\mathfrak{so}(14)_{\mathbb{C}}}$,
$[X_{120}]_{\mathfrak{so}(14)_{\mathbb{C}}}$, $[X_{112}]_{\mathfrak{so}(14)_{\mathbb{C}}}$,
$[Y_{97}]_{\mathfrak{so}(14)_{\mathbb{C}}}$  and~$[Y_1]_{\mathfrak{so}(14)_{\mathbb{C}}}$ from equation~\eqref{hhhhkj}.

}

Let $\alpha_1, \alpha_2, \alpha_3,\ldots,\alpha_8$
be a~set of simple roots for~$E_8$.
To any positive root $a_1 \alpha_1+a_2\alpha_2 +a_3 \alpha_3+ \dots +a_8\alpha_8$
we may associate a~vector $[a_1, a_2, a_3,\ldots,a_8]\in \mathbb{Z}^{8}_{\geq 0}$.
With this convention, the positive roots of $E_8$, as computed with GAP~\cite{gap}, are as follows:
\begin{alignat*}{3}
& \alpha_{1} = [1, 0, 0, 0, 0, 0, 0, 0], \qquad && \alpha_{2} =  [0, 1, 0, 0, 0, 0, 0, 0],&
\\
&\alpha_{3} =  [0, 0, 1, 0, 0, 0, 0, 0], \qquad && \alpha_{4} =  [0, 0, 0, 1, 0, 0, 0, 0],&
\\
&\alpha_{5} =  [0, 0, 0, 0, 1, 0, 0, 0], \qquad && \alpha_{6} =  [0, 0, 0, 0, 0, 1, 0, 0],&
\\
&\alpha_{7} =  [0, 0, 0, 0, 0, 0, 1, 0], \qquad && \alpha_{8} =  [0, 0, 0, 0, 0, 0, 0, 1],&
\\
&\alpha_{9} =  [1, 0, 1, 0, 0, 0, 0, 0], \qquad && \alpha_{10} =  [0, 1, 0, 1, 0, 0, 0, 0],&
\\
&\alpha_{11} =  [0, 0, 1, 1, 0, 0, 0, 0], \qquad && \alpha_{12} =  [0, 0, 0, 1, 1, 0, 0, 0],&
\\
& \alpha_{13} =  [0, 0, 0, 0, 1, 1, 0, 0], \qquad && \alpha_{14} =  [0, 0, 0, 0, 0, 1, 1, 0],&
\\
&\alpha_{15} =  [0, 0, 0, 0, 0, 0, 1, 1], \qquad && \alpha_{16} =  [1, 0, 1, 1, 0, 0, 0, 0],&
\\
&\alpha_{17} =  [0, 1, 1, 1, 0, 0, 0, 0], \qquad && \alpha_{18} =  [0, 1, 0, 1, 1, 0, 0, 0],&
\\
&\alpha_{19} =  [0, 0, 1, 1, 1, 0, 0, 0], \qquad && \alpha_{20} =  [0, 0, 0, 1, 1, 1, 0, 0],&
\\
&\alpha_{21} =  [0, 0, 0, 0, 1, 1, 1, 0], \qquad && \alpha_{22} =  [0, 0, 0, 0, 0, 1, 1, 1],&
\\
&\alpha_{23} =  [1, 1, 1, 1, 0, 0, 0, 0], \qquad && \alpha_{24} =  [1, 0, 1, 1, 1, 0, 0, 0],&
\\
&\alpha_{25} =  [0, 1, 1, 1, 1, 0, 0, 0], \qquad && \alpha_{26} =  [0, 1, 0, 1, 1, 1, 0, 0],&
\\
&\alpha_{27} =  [0, 0, 1, 1, 1, 1, 0, 0], \qquad && \alpha_{28} =  [0, 0, 0, 1, 1, 1, 1, 0],&
\\
&\alpha_{29} =  [0, 0, 0, 0, 1, 1, 1, 1], \qquad && \alpha_{30} =  [1, 1, 1, 1, 1, 0, 0, 0],&
\\
&\alpha_{31} =  [1, 0, 1, 1, 1, 1, 0, 0], \qquad && \alpha_{32} =  [0, 1, 1, 2, 1, 0, 0, 0],&
\\
&\alpha_{33} =  [0, 1, 1, 1, 1, 1, 0, 0], \qquad && \alpha_{34} =  [0, 1, 0, 1, 1, 1, 1, 0],&
\\
&\alpha_{35} =  [0, 0, 1, 1, 1, 1, 1, 0], \qquad && \alpha_{36} =  [0, 0, 0, 1, 1, 1, 1, 1],&
\\
&\alpha_{37} =  [1, 1, 1, 2, 1, 0, 0, 0], \qquad && \alpha_{38} =  [1, 1, 1, 1, 1, 1, 0, 0],&
\\
&\alpha_{39} =  [1, 0, 1, 1, 1, 1, 1, 0], \qquad && \alpha_{40} =  [0, 1, 1, 2, 1, 1, 0, 0],&
\\
&\alpha_{41} =  [0, 1, 1, 1, 1, 1, 1, 0], \qquad && \alpha_{42} =  [0, 1, 0, 1, 1, 1, 1, 1],&
\\
&\alpha_{43} =  [0, 0, 1, 1, 1, 1, 1, 1], \qquad && \alpha_{44} =  [1, 1, 2, 2, 1, 0, 0, 0],&
\\
&\alpha_{45} =  [1, 1, 1, 2, 1, 1, 0, 0], \qquad && \alpha_{46} =  [1, 1, 1, 1, 1, 1, 1, 0],&
\\
&\alpha_{47} =  [1, 0, 1, 1, 1, 1, 1, 1], \qquad && \alpha_{48} =  [0, 1, 1, 2, 2, 1, 0, 0],&
\\
&\alpha_{49} =  [0, 1, 1, 2, 1, 1, 1, 0], \qquad && \alpha_{50} =  [0, 1, 1, 1, 1, 1, 1, 1],&
\\
&\alpha_{51} =  [1, 1, 2, 2, 1, 1, 0, 0], \qquad && \alpha_{52} =  [1, 1, 1, 2, 2, 1, 0, 0],&
\\
&\alpha_{53} =  [1, 1, 1, 2, 1, 1, 1, 0], \qquad && \alpha_{54} =  [1, 1, 1, 1, 1, 1, 1, 1],&
\\
&\alpha_{55} =  [0, 1, 1, 2, 2, 1, 1, 0], \qquad && \alpha_{56} =  [0, 1, 1, 2, 1, 1, 1, 1],&
\\
&\alpha_{57} =  [1, 1, 2, 2, 2, 1, 0, 0], \qquad && \alpha_{58} =  [1, 1, 2, 2, 1, 1, 1, 0],&
\\
&\alpha_{59} =  [1, 1, 1, 2, 2, 1, 1, 0], \qquad && \alpha_{60} =  [1, 1, 1, 2, 1, 1, 1, 1],&
\\
&\alpha_{61} =  [0, 1, 1, 2, 2, 2, 1, 0], \qquad && \alpha_{62} =  [0, 1, 1, 2, 2, 1, 1, 1],&
\\
&\alpha_{63} =  [1, 1, 2, 3, 2, 1, 0, 0], \qquad && \alpha_{64} =  [1, 1, 2, 2, 2, 1, 1, 0],&
\\
&\alpha_{65} =  [1, 1, 2, 2, 1, 1, 1, 1], \qquad && \alpha_{66} =  [1, 1, 1, 2, 2, 2, 1, 0],&
\\
&\alpha_{67} =  [1, 1, 1, 2, 2, 1, 1, 1], \qquad && \alpha_{68} =  [0, 1, 1, 2, 2, 2, 1, 1],&
\\
&\alpha_{69} =  [1, 2, 2, 3, 2, 1, 0, 0], \qquad && \alpha_{70} =  [1, 1, 2, 3, 2, 1, 1, 0],&
\\
&\alpha_{71} =  [1, 1, 2, 2, 2, 2, 1, 0], \qquad && \alpha_{72} =  [1, 1, 2, 2, 2, 1, 1, 1],&
\\
&\alpha_{73} =  [1, 1, 1, 2, 2, 2, 1, 1], \qquad && \alpha_{74} =  [0, 1, 1, 2, 2, 2, 2, 1],&
\\
&\alpha_{75} =  [1, 2, 2, 3, 2, 1, 1, 0], \qquad && \alpha_{76} =  [1, 1, 2, 3, 2, 2, 1, 0],&
\\
&\alpha_{77} =  [1, 1, 2, 3, 2, 1, 1, 1], \qquad && \alpha_{78} =  [1, 1, 2, 2, 2, 2, 1, 1],&
\\
&\alpha_{79} =  [1, 1, 1, 2, 2, 2, 2, 1], \qquad && \alpha_{80} =  [1, 2, 2, 3, 2, 2, 1, 0],&
\\
&\alpha_{81} =  [1, 2, 2, 3, 2, 1, 1, 1], \qquad && \alpha_{82} =  [1, 1, 2, 3, 3, 2, 1, 0],&
\\
&\alpha_{83} =  [1, 1, 2, 3, 2, 2, 1, 1], \qquad && \alpha_{84} =  [1, 1, 2, 2, 2, 2, 2, 1],&
\\
&\alpha_{85} =  [1, 2, 2, 3, 3, 2, 1, 0], \qquad && \alpha_{86} =  [1, 2, 2, 3, 2, 2, 1, 1],&
\\
&\alpha_{87} =  [1, 1, 2, 3, 3, 2, 1, 1], \qquad && \alpha_{88} =  [1, 1, 2, 3, 2, 2, 2, 1],&
\\
&\alpha_{89} =  [1, 2, 2, 4, 3, 2, 1, 0], \qquad && \alpha_{90} =  [1, 2, 2, 3, 3, 2, 1, 1],&
\\
&\alpha_{91} =  [1, 2, 2, 3, 2, 2, 2, 1], \qquad && \alpha_{92} =  [1, 1, 2, 3, 3, 2, 2, 1],&
\\
&\alpha_{93} =  [1, 2, 3, 4, 3, 2, 1, 0], \qquad && \alpha_{94} =  [1, 2, 2, 4, 3, 2, 1, 1],&
\\
&\alpha_{95} =  [1, 2, 2, 3, 3, 2, 2, 1], \qquad && \alpha_{96} =  [1, 1, 2, 3, 3, 3, 2, 1],&
\\
&\alpha_{97} =  [2, 2, 3, 4, 3, 2, 1, 0], \qquad && \alpha_{98} =  [1, 2, 3, 4, 3, 2, 1, 1],&
\\
&\alpha_{99} =  [1, 2, 2, 4, 3, 2, 2, 1], \qquad && \alpha_{100} =  [1, 2, 2, 3, 3, 3, 2, 1],&
\\
&\alpha_{101} =  [2, 2, 3, 4, 3, 2, 1, 1], \qquad && \alpha_{102} =  [1, 2, 3, 4, 3, 2, 2, 1],&
\\
&\alpha_{103} =  [1, 2, 2, 4, 3, 3, 2, 1], \qquad && \alpha_{104} =  [2, 2, 3, 4, 3, 2, 2, 1],&
\\
&\alpha_{105} =  [1, 2, 3, 4, 3, 3, 2, 1], \qquad && \alpha_{106} =  [1, 2, 2, 4, 4, 3, 2, 1],&
\\
&\alpha_{107} =  [2, 2, 3, 4, 3, 3, 2, 1], \qquad && \alpha_{108} =  [1, 2, 3, 4, 4, 3, 2, 1],&
\\
&\alpha_{109} =  [2, 2, 3, 4, 4, 3, 2, 1], \qquad && \alpha_{110} =  [1, 2, 3, 5, 4, 3, 2, 1],&
\\
&\alpha_{111} =  [2, 2, 3, 5, 4, 3, 2, 1], \qquad && \alpha_{112} =  [1, 3, 3, 5, 4, 3, 2, 1],&
\\
&\alpha_{113} =  [2, 3, 3, 5, 4, 3, 2, 1], \qquad && \alpha_{114} =  [2, 2, 4, 5, 4, 3, 2, 1],&
\\
&\alpha_{115} =  [2, 3, 4, 5, 4, 3, 2, 1], \qquad && \alpha_{116} =  [2, 3, 4, 6, 4, 3, 2, 1],&
\\
&\alpha_{117} =  [2, 3, 4, 6, 5, 3, 2, 1], \qquad && \alpha_{118} =  [2, 3, 4, 6, 5, 4, 2, 1],&
\\
&\alpha_{119} =  [2, 3, 4, 6, 5, 4, 3, 1], \qquad && \alpha_{120} =  [2, 3, 4, 6, 5, 4, 3, 2].&
\end{alignat*}

Let $X_{\alpha_i}=X_i$, and $Y_{\alpha_i}=Y_i$ be a~choice of positive (resp.
negative) root vector corresponding to the root $\alpha_i$.

A basis of $V(\lambda_1)=[X_{120}]_{\mathfrak{so}(14)_{\mathbb{C}}}$ is given by the $14$ positive root vectors:
\begin{alignat*}{8}
& X_{97}, \quad&& X_{101}, \quad&& X_{104}, \quad&& X_{107}, \quad&& X_{109}, \quad&& X_{111}, \quad&& X_{113},&
\\
& X_{114}, \quad&& X_{115}, \quad&& X_{116}, \quad&& X_{117}, \quad&& X_{118}, \quad&& X_{119}, \quad&& X_{120}.&
\end{alignat*}

A basis of $V(\lambda_1)=[Y_{97}]_{\mathfrak{so}(14)_{\mathbb{C}}}$ is given by the $14$ negative root vectors:
\begin{alignat*}{8}
& Y_{97}, \quad&& Y_{101}, \quad&& Y_{104}, \quad&& Y_{107}, \quad&& Y_{109}, \quad&& Y_{111}, \quad&& Y_{113},&
\\
& Y_{114}, \quad&& Y_{115}, \quad&& Y_{116}, \quad&& Y_{117}, \quad&& Y_{118}, \quad&& Y_{119}, \quad&& Y_{120}.&
\end{alignat*}

A basis of $V(\lambda_7)=[X_{112}]_{\mathfrak{so}(14)_{\mathbb{C}}}$ is given by the $64$ positive root vectors:
\begin{alignat*}{10}
&X_{1}, \quad&& X_{9}, \quad&& X_{16}, \quad&& X_{23}, \quad&& X_{24}, \quad&& X_{30}, \quad&& X_{31}, \quad&& X_{37}, \quad&& X_{38},&
\\
&X_{39}, \quad&& X_{44}, \quad&& X_{45}, \quad&& X_{46}, \quad&& X_{47}, \quad&& X_{51}, \quad&& X_{52}, \quad&& X_{53}, \quad&& X_{54},&
\\
&X_{57}, \quad&& X_{58}, \quad&& X_{59}, \quad&& X_{60}, \quad&& X_{63}, \quad&& X_{64}, \quad&& X_{65}, \quad&& X_{66}, \quad&& X_{67},&
\\
&X_{69}, \quad&& X_{70}, \quad&& X_{71}, \quad&& X_{72}, \quad&& X_{73}, \quad&& X_{75}, \quad&& X_{76}, \quad&& X_{77}, \quad&& X_{78},&
\\
&X_{79}, \quad&& X_{80}, \quad&& X_{81}, \quad&& X_{82}, \quad&& X_{83}, \quad&& X_{84}, \quad&& X_{85}, \quad&& X_{86}, \quad&& X_{87},&
\\
&X_{88}, \quad&& X_{89}, \quad&& X_{90}, \quad&& X_{91}, \quad&& X_{92}, \quad&& X_{93}, \quad&& X_{94}, \quad&& X_{95}, \quad&& X_{96},&
\\
&X_{98}, \quad&& X_{99}, \quad&& X_{100}, \quad&& X_{102}, \quad&& X_{103}, \quad&& X_{105}, \quad&& X_{106}, \quad&& X_{108}, \quad&& X_{110},&
\\
&X_{112}. &&&&&&&&&&&&&&&& &
\end{alignat*}

A basis of $V(\lambda_6)=[Y_1]_{\mathfrak{so}(14)_{\mathbb{C}}}$ is given by the $64$ negative root vectors:
\begin{alignat*}{10}
&Y_{1}, \quad&& Y_{9}, \quad&& Y_{16}, \quad&& Y_{23}, \quad&& Y_{24}, \quad&& Y_{30}, \quad&& Y_{31}, \quad&& Y_{37}, \quad&& Y_{38},&
\\
&Y_{39}, \quad&& Y_{44}, \quad&& Y_{45}, \quad&& Y_{46}, \quad&& Y_{47}, \quad&& Y_{51}, \quad&& Y_{52}, \quad&& Y_{53}, \quad&& Y_{54},&
\\
& Y_{57}, \quad&& Y_{58}, \quad&& Y_{59}, \quad&& Y_{60}, \quad&& Y_{63}, \quad&& Y_{64}, \quad&& Y_{65}, \quad&& Y_{66}, \quad&& Y_{67},&
\\
&Y_{69}, \quad&& Y_{70}, \quad&& Y_{71}, \quad&& Y_{72}, \quad&& Y_{73}, \quad&& Y_{75}, \quad&& Y_{76}, \quad&& Y_{77}, \quad&& Y_{78},&
\\
&Y_{79}, \quad&& Y_{80}, \quad&& Y_{81}, \quad&& Y_{82}, \quad&& Y_{83}, \quad&& Y_{84}, \quad&& Y_{85}, \quad&& Y_{86}, \quad&& Y_{87},&
\\
&Y_{88}, \quad&& Y_{89}, \quad&& Y_{90}, \quad&& Y_{91}, \quad&& Y_{92}, \quad&& Y_{93}, \quad&& Y_{94}, \quad&& Y_{95}, \quad&& Y_{96},&
\\
&Y_{98}, \quad&& Y_{99}, \quad&& Y_{100}, \quad&& Y_{102}, \quad&& Y_{103}, \quad&& Y_{105}, \quad&& Y_{106}, \quad&& Y_{108}, \quad&& Y_{110},&
\\
&Y_{112}. &&&&&&&&&&&&&&&& &
\end{alignat*}

Note that $\varphi(\mathfrak{so}(14)_{\mathbb{C}}) =[X_{74}]_{\mathfrak{so}(14)_{\mathbb{C}}}$.
We describe bases of $\varphi({\mathfrak{so}(14)_{\mathbb{C}+}})$ and $\varphi({\mathfrak{so}(14)_{\mathbb{C}-}})$,
respectively:
\begin{alignat*}{10}
&X_{2}, \quad&&  X_{3}, \quad&&  X_{4}, \quad&&  X_{5}, \quad&&  X_{6}, \quad&&  X_{7}, \quad&&  X_{8}, \quad&&  X_{10}, \quad&&  X_{11},&
\\
&X_{12}, \quad&&  X_{13}, \quad&&  X_{14}, \quad&&  X_{15}, \quad&&  X_{17}, \quad&&  X_{18}, \quad&&  X_{19}, \quad&&  X_{20}, \quad&&  X_{21},&
\\
&X_{22}, \quad&&  X_{25}, \quad&&  X_{26}, \quad&&  X_{27}, \quad&&  X_{28}, \quad&&  X_{29}, \quad&&  X_{32}, \quad&&  X_{33}, \quad&&  X_{34},&
\\
&X_{35}, \quad&&  X_{36}, \quad&&  X_{40}, \quad&&  X_{41}, \quad&&  X_{42}, \quad&&  X_{43}, \quad&&  X_{48}, \quad&&  X_{49}, \quad&&  X_{50},&
\\
&X_{55}, \quad&&  X_{56}, \quad&&  X_{61}, \quad&&  X_{62}, \quad&&  X_{68}, \quad&&  X_{74}; && && && &
\\
&Y_{2}, \quad&&  Y_{3}, \quad&&  Y_{4}, \quad&&  Y_{5}, \quad&&  Y_{6}, \quad&&  Y_{7}, \quad&&  Y_{8}, \quad&&  Y_{10}, \quad&&  Y_{11},&
\\
&Y_{12}, \quad&&  Y_{13}, \quad&&  Y_{14}, \quad&&  Y_{15}, \quad&&  Y_{17}, \quad&&  Y_{18}, \quad&&  Y_{19}, \quad&&  Y_{20}, \quad&&  Y_{21},&
\\
&Y_{22}, \quad&&  Y_{25}, \quad&&  Y_{26}, \quad&&  Y_{27}, \quad&&  Y_{28}, \quad&&  Y_{29}, \quad&&  Y_{32}, \quad&&  Y_{33}, \quad&&  Y_{34},&
\\
&Y_{35}, \quad&&  Y_{36}, \quad&&  Y_{40}, \quad&&  Y_{41}, \quad&&  Y_{42}, \quad&&  Y_{43}, \quad&&  Y_{48}, \quad&&  Y_{49}, \quad&&  Y_{50},&
\\
&Y_{55}, \quad&&  Y_{56}, \quad&&  Y_{61}, \quad&&  Y_{62}, \quad&&  Y_{68}, \quad&&  Y_{74}. &&&&&& &
\end{alignat*}

{\bf Acknowledgements.}
The work of A.D.\
is partially supported by a~research grant from the Professional Staf\/f Congress/ City University of New York (PSC/CUNY).
The work of J.R.\
is partially supported by the Natural Sciences and Engineering Research Council (NSERC).
The authors would also like to thank the anonymous referees for valuable comments.

\vspace{-2mm}

\pdfbookmark[1]{References}{ref}
\LastPageEnding

\end{document}